\documentclass[12pt]{article}
\usepackage{amsmath,amsthm,amssymb,amsfonts}
\usepackage[numbers,sort&compress]{natbib}
\usepackage{color,colordvi}
\usepackage{soul}
\usepackage{fullpage}

\usepackage{amsmath,amsfonts,amsthm,mathtools,bm}
\usepackage{graphicx}
\usepackage[colorlinks=true, allcolors=blue]{hyperref}
\usepackage[capitalise]{cleveref}
\usepackage{standalone}
\usepackage{indentfirst}
\usepackage[affil-it]{authblk}
\usepackage{diagbox}
\usepackage{changes}
\usepackage{enumitem}
\usepackage{bbm}
\usepackage{orcidlink}

\newtheorem{theorem}{Theorem}[section]
\newtheorem{lemma}[theorem]{Lemma}
\newtheorem{corollary}[theorem]{Corollary}

\newtheorem{conjecture}[theorem]{Conjecture}

\newtheorem{definition}[theorem]{Definition}

\newtheorem*{claim*}{Claim}

\newcommand{\textotherwise}{\text{otherwise}}

\newcommand{\RR}{\mathbb{R}}

\newcommand{\ZZ}{\mathbb{Z}}

\newcommand{\cT}{\mathcal{T}}

\DeclarePairedDelimiter{\card}{\lvert}{\rvert}
\DeclarePairedDelimiter{\set}{\lbrace}{\rbrace}

\title{Maximize the Steklov eigenvalue of trees}

\author{Huiqiu Lin\footnote{email: huiqiulin@126.com}}
\author{Da Zhao\footnote{email: zhaoda@ecust.edu.cn}~\orcidlink{0000-0002-9582-0778}}
\affil{School of Mathematics, East China University of Science and Technology, 130 Meilong Road, Shanghai 200237, China.}

\date{}

\begin{document}
\maketitle

\begin{abstract} 
    We study the maximal Steklov eigenvalues of trees with given number of boundary vertices and total number of vertices. 
    Trees can be regarded as discrete analogue of Hadamard manifolds, namely simply-connected Riemannian manifolds of non-positive sectional curvature. 
    Let $\sigma_{k,\text{max}}(b, n)$ be the maximal of $k$-th Steklov eigenvalue of trees with $b$ leaves as boundary and $n$ vertices. 
    We determine that 
    $$
    \sigma_{2, \text{max}} (b, n) = \begin{cases}
        \frac{2}{n-1}, & b=2, n\geq 3, \\
        \frac{1}{r}, & b \geq 3, n = br + m, 3 - b \leq m \leq 1, r \in \mathbb{Z}_+, \\
        \frac{1}{r+1-\frac{1}{b}}, & b \geq 3, n = br + 2, r \in \mathbb{Z}_+,
    \end{cases}
    $$
    and we characterize the trees attaining this bound. 
    For $k \geq 3$, we show that $\sigma_{k, \text{max}} (b, n) = 1$. 
    We also give a lower bound on the maximal Steklov eigenvalues of trees with given diameter and total number of vertices. 
    Our work can be regarded as a completion of the work by He--Hua [Upper bounds
    for the Steklov eigenvalues on trees, Calc. Var. Partial Differential Equations (2022)] and Yu--Yu [Monotonicity of Steklov eigenvalues on graphs
    and applications, Calc. Var. Partial Differential Equations (2024)].
\end{abstract}


\section{Introduction}

Steklov~\cite{stekloff1902ProblemesFondamentauxPhysique, kuznetsov2014LegacyVladimirAndreevich} introduced the Steklov problem over a century ago during the study of liquid sloshing. 
Given a compact Riemannian manifold with boundary, the Steklov operator is defined to be the map sending the Dirichlet boundary data of a harmonic function on the manifold to its Neumann boundary data, thus it is also called Dirichlet-to-Neumann operator, abbreviated by DtN operator. 
The Steklov operator establishes the model for inverse problem of detecting the inside of a body by measuring the data on the boundary. 
The eigenvalues of the Steklov operators are called Steklov eigenvalues. 
The relationship between Steklov eigenvalues and the Yamabe problem on Riemannian manifold with boundary was found by Escobar~\cite{escobar1992YamabeProblemManifolds}. 
In the last decade, Fraser and Scheon~\cite{fraser2011FirstSteklovEigenvalue, fraser2016SharpEigenvalueBounds, fraser2019ShapeOptimizationSteklov, fraser2020ResultsHigherEigenvalue} showed the relationship between Steklov eigenfunctions and free boundary minimal submanifolds in Euclidean ball and proceed on extremal problems of Steklov eigenvalues. 

Recently Hua--Huang--Wang~\cite{hua2017FirstEigenvalueEstimates} and Hassannezhad--Miclo~\cite{asmahassannezhad2020HigherOrderCheeger} independently extended the Steklov problem to discrete spaces. 
Colbois--Girouard~\cite{girouard2014SpectralGapGraphs} on the other hand constructed compact surfaces with first nontrivial Steklov eigenvalue uniformly bounded away from zero from expander graphs. 
Later numerous studies investigated the properties of Steklov eigenvalues on graphs. 
In~\cite{han2023SteklovEigenvalueProblem, asmahassannezhad2020HigherOrderCheeger,hua2017FirstEigenvalueEstimates, hua2023CheegerEstimatesDirichlettoneumann, perrin2019LowerBoundsFirst, perrin2020IsoperimetricUpperBound,tschanz2022UpperBoundsSteklov}, the authors controlled the Steklov eigenvalues by isoperimetrical parameters. 
He--Hua studied the upper bound of first Steklov eigenvalue on trees in~\cite{he2022UpperBoundsSteklov,he2022SteklovFlowsTrees}. 
In~\cite{yu2024MonotonicitySteklovEigenvalues}, Yu--Yu proved the monotonicity of Steklov eigenvalues. 
Perrin provided some lower bounds of Steklov eigenvalues in~\cite{perrin2019LowerBoundsFirst}, and Shi--Yu generalized the bound to weighted graphs~\cite{shi2025ExtensionRigidityPerrins}. 
Shi--Yu compared the Steklov eigenvalue and the Laplacian eigenvalue in~\cite{shi2022ComparisonSteklovEigenvalues}. 
In~\cite{shi2022LichnerowicztypeEstimateSteklov}, the authors proved a Lichnerowicz--type estimate for Steklov eigenvalues. 
In~\cite{shi2022DirichlettoneumannMapsDifferential}, the authors established Steklov operators for differential forms on graphs. 
The survey paper \cite{colbois2023RecentDevelopmentsSteklov} summarizes recent progress on Steklov problem. 

In~\cite{yu2022MinimalSteklovEigenvalues}, the authors characterized the minimal Steklov eigenvalues on graphs. 
In this paper, we characterized the maximal Steklov eigenvalues on trees, which can be regarded as a completion of the work in~\cite{he2022UpperBoundsSteklov,he2022SteklovFlowsTrees,yu2024MonotonicitySteklovEigenvalues}.

A graph is a pair $G = (V, E)$, where $V$ is the vertex set and $E \subseteq \binom{V}{2}$ is the edge set.
The edge $(x, y) \in E$ is sometimes abbreviated as $xy$ or $x \sim y$. 
Let $B$ be a subset of the vertex set, which we call the boundary of the graph. 

For a graph $G$ with boundary $B$, we can define the (discrete) Steklov operator $\Lambda: \RR^B \to \RR^B$ by $\Lambda(f) = \frac{\partial \hat{f}}{\partial n}$, where $\hat{f}$ is the harmonic extension of $f$ to the whole graph $V$, and $\frac{\partial \hat{f}}{\partial n}$ is the discrete normal derivative given by $\frac{\partial \hat{f}}{\partial n}(x) = \sum_{(x, y) \in E} (f(x) - f(y))$. 
The eigenvalues of $\Lambda$ are called Steklov eigenvalues of the pair $(G, B)$, denoted by 
\begin{align}
    0 = \sigma_1(G,B) \leq \sigma_2(G,B) \leq \cdots \leq \sigma_{\card{B}}(G, B). 
\end{align}
We write $\sigma_i$ for $\sigma_i(G, B)$ when the graph and the boundary is clear from context. 

A path of length $\ell$ in a graph is a sequence of distinct vertices $v_0 v_1 \cdots v_{\ell}$ such that $v_{i-1} \sim v_i$ for all $i = 1,2, \ldots, \ell$. 
The distance between two vertices $x, y \in V$ is the length of shorted path connecting them. 
The diameter $D$ of a graph $G$ is the maximum distance among all pairs of vertices. 
A graph $G$ is connected if for every pair of vertices $x, y \in V$, there exists a path connecting them. 
A tree is a minimal connected graph, namely it is a connected graph with $(\card{V} - 1)$ edges. 
A path graph is a graph consisting of a path. 
The degree of a vertex $v \in V$ is the number of edges adjacent to it, namely $\deg(v) = \card{\set{u \in V: u \sim v}}$. 
The vertices of degree $1$ in a tree are called leaves. 
We take the set of leaves as the natural boundary of a tree. 

Let $k, b, n$ be integers such that $1 \leq k \leq b$ and $n \geq b+1$. 
Consider the parameter $\sigma_{k, \text{max}}(b, n)$ defined by the maximum $k$-th Steklov eigenvalue among all trees with $b$ leaves and $n$ vertices, namely $\sigma_{k, \text{max}}(b, n) = \max_{T \in \cT(b, n)} \sigma_k(T)$, where $\cT(b, n)$ is the set of trees with $b$ leaves and $n$ vertices. 
Note that $\sigma_{1, \text{max}}(b, n) = 0$. 

Our first main result is the exact value of $\sigma_{2, \text{max}}(b, n)$.

\begin{theorem}\label{thm:question}
    Let $b \geq 2$ and $n \geq b+1$ be integers. 
    Then the followings hold. 
    \begin{align}\label{eq:parent}
        \sigma_{2, \text{max}}(2, n) = \frac{2}{n-1},
    \end{align}
    and for $b \geq 3$, it holds
    \begin{align}
        \sigma_{2, \text{max}}(b, n) = 
        \begin{cases}
            \frac{1}{r}, & n = br + m, 3-b \leq m \leq 1, r \in \ZZ_+, \\
            \frac{1}{r + 1 - \frac{1}{b}}, & n = br+2, r \in \ZZ_+.
        \end{cases}
    \end{align}

    The tree attaining the bound $\sigma_{2, \text{max}} = \frac{1}{r}$  is not unique. 
    Let $v_0, v_1, \ldots, v_{2r} = v_{D}$ be a diametral path of length $2r$. 
    For $i = 1,2, \ldots, D-1$, let $H_i$ be the connected components containing neither $v_0$ nor $v_{D}$ obtained by deleting $v_{i-1}v_i$ and $v_i v_{i+1}$ from edges of $T$. 
    As long as $H_i$ is empty for all $1 \leq i \leq D-1$ except for $i = r$, and $H_r$ is a rooted tree with depth at most $r$, the first (nontrivial) Steklov eigenvalue of the tree $T$ is $\frac{1}{r}$. 
    The tree attaining the bound $\sigma_{2, \text{max}} = \frac{1}{r+1 - \frac{1}{b}}$ is unique. 
    It is exactly the almost fork graph $AF(b, r)$ (see~\cref{def:nearby}).
\end{theorem}

\begin{corollary}
    Let $T = (V, E)$ be a tree with leaves as boundary $B$. 
    Let $D$ be the diameter of $T$. 
    Then 
    \begin{align}
        \sigma_2 \leq \frac{2}{D}.
    \end{align}
    Moreover, if $D$ is odd, then a necessary and sufficient condition to attain the upper bound is that $T$ is a path graph. 
    For $i = 1,2, \ldots, D-1$, let $H_i$ be the connected components containing neither $v_0$ nor $v_{D}$ obtained by deleting $v_{i-1}v_i$ and $v_i v_{i+1}$ from edges of $T$. 
    If $D = 2r$ is even, then a necessary and sufficient condition to attain the upper bound is that $H_i$ is empty for all $1 \leq i \leq D-1$ except for $i = r$, and $H_r$ is a rooted tree of depth at most $r$. 
\end{corollary}

We also determine the exact value for $\sigma_{k, \text{max}}(b, n)$ for $3 \leq k \leq b$.

\begin{theorem}\label{thm:natural}
    Let $b \geq 2$ and $n \geq b+1$ be integers. 
    Then $\sigma_{k, \text{max}}(b, n) = 1$ for $3 \leq k \leq b$. 
\end{theorem}

Let $D, n$ be integers such that $n \geq D+1 \geq 2$. 
Consider the parameter $\widetilde{\sigma}_{2, \text{max}}(D, n)$ defined by the maximum first nontrivial Steklov eigenvalue among all trees with diameter $D$ and $n$ vertices, namely $\widetilde{\sigma}_{2, \text{max}}(D, n) = \max_{T \in \widetilde{\cT}(D, n)} \sigma_2(T)$, where $\widetilde{\cT}(D, n)$ is the set of trees with diameter $D$ and $n$ vertices. 

We obtain lower bounds for the value of $\widetilde{\sigma}_{2, \text{max}}(D, n)$. 

\begin{theorem}\label{thm:gas}
    Let $D \geq 1$ and $n \geq D+1$ be integers. 
    Then the followings hold. 
    \begin{align}\label{eq:order}
        \widetilde{\sigma}_{2, \text{max}}(1, 2) = 1,
    \end{align}
    and for $D \geq 2$, it holds
    \begin{align}
        \widetilde{\sigma}_{2, \text{max}}(D, n) 
        \begin{cases}
            = \frac{2}{D}, & D=2r \text{ is even}, \\
            = \frac{n-2}{2n-5}, & D=3, \\
            \geq C^-(r, b, c)& D=2r+1 \geq 5 \text{ is odd, } \\
            & c \leq r, \text{ and } 2r + 2 \leq n \leq 2r + (b-2)c + 2. 
        \end{cases}
    \end{align}
    The definition of $C^-(r, b, c)$ is given in~\cref{lem:almost_seesaw}.
    The tree attaining the bound $\widetilde{\sigma}_{2, \text{max}} = \frac{2}{D}$ is not unique. 
    As long as $H_r$ is a rooted tree with depth at most $r = \frac{1}{2}D$, the first Steklov eigenvalue of the tree $T$ is $\frac{2}{D}$. 
\end{theorem}

\section{Proofs}

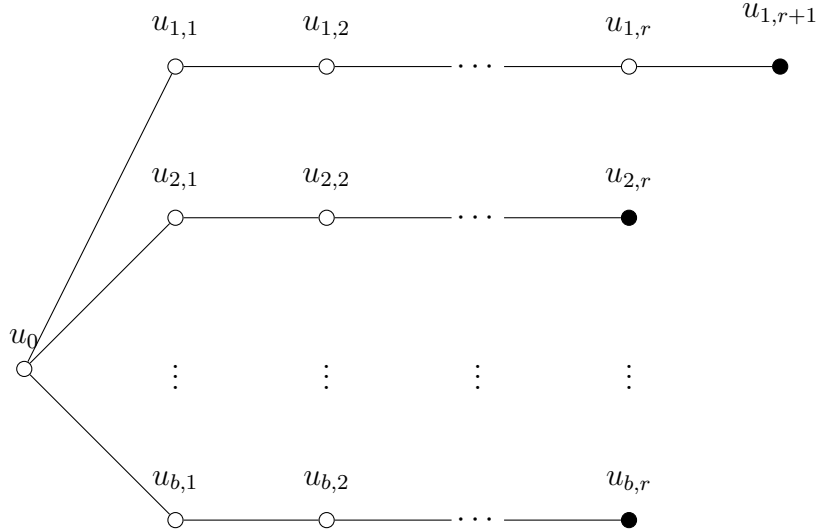
\begin{figure}[htbp]
    \centering
    \begin{tikzpicture}
    [scale=2, every node/.style={circle, inner sep=1pt, minimum size = 2mm}]

    \node[draw, label=above:{$u_0$}] (u0) at (0, 0) {};
    \node[draw, label=above:{$u_{1,1}$}] (u11) at (1, 2) {};
    \node[draw, label=above:{$u_{2,1}$}] (u21) at (1, 1) {};
    \node (um1) at (1, 0) {$\vdots$};
    \node[draw, label=above:{$u_{b,1}$}] (ub1) at (1, -1) {};

    \node[draw, label=above:{$u_{1,2}$}] (u12) at (2, 2) {};
    \node[draw, label=above:{$u_{2,2}$}] (u22) at (2, 1) {};
    \node (um2) at (2, 0) {$\vdots$};
    \node[draw, label=above:{$u_{b,2}$}] (ub2) at (2, -1) {};

    \node (u13) at (3, 2) {$\cdots$};
    \node (u23) at (3, 1) {$\cdots$};
    \node (um3) at (3, 0) {$\vdots$};
    \node (ub3) at (3, -1) {$\cdots$};

    \node[draw, label=above:{$u_{1,r}$}] (u14) at (4, 2) {};
    \node[draw, fill=black, label=above:{$u_{2,r}$}] (u24) at (4, 1) {};
    \node (um4) at (4, 0) {$\vdots$};
    \node[draw, fill=black, label=above:{$u_{b,r}$}] (ub4) at (4, -1) {};

    \node[draw, fill=black, label=above:{$u_{1,r+1}$}] (u15) at (5, 2) {};

    \draw (u0) -- (u11) -- (u12) -- (u13) -- (u14) -- (u15);
    \draw (u0) -- (u21) -- (u22) -- (u23) -- (u24);
    \draw (u0) -- (ub1) -- (ub2) -- (ub3) -- (ub4);

\end{tikzpicture}
    \caption{an almost fork graph $AF(b, r)$ with leaves as boundary.}\label{fig:luck}
\end{figure}

\begin{definition}\label{def:nearby}
    Let $b \geq 2$ and $r \geq 1$ be integers. 
    The almost fork graph $AF(b, r)$ is obtained by gluing one ends of $b$ paths together, where one path is of length $r+1$ and the rest $(b-1)$ paths are of length $r$. 
    See~\cref{fig:luck}.
\end{definition}

We give the Steklov eigenvalues of almost fork graphs. 

\begin{lemma}
    Let $T = (V, E)$ be the almost fork graph $AF(b, r)$ with leaves as boundary $B$. 
    Then the Steklov eigenvalues of $G$ are
    $\sigma_1 = 0$, $\sigma_2 = \frac{b}{b(r+1)-1}$, and $\sigma_3 = \cdots = \sigma_b = 1/r$ for $b \geq 2$ and $r \geq 1$. 
\end{lemma}

\begin{proof}
    We give the Steklov eigenfunctions directly. 
    Define $\sigma_1 = 0$, $\sigma_2 = \frac{b}{b(r+1)-1}$, and $\sigma_3 = \cdots = \sigma_b = 1/r$. 
    Define
    \begin{align}
        \xi_1(v) = 1, \quad \forall v \in V; 
    \end{align}
    \begin{align}
        \xi_2(v) = 
        \begin{cases}
            -\frac{b-1}{b}\sigma_2, & v = u_0, \\
            -(b-1)(1 - (r+1 - i)\sigma_2), & v = u_{1, i},\ 1 \leq i \leq r+1, \\
            1 - (r-i)\sigma_2, & v = u_{\ell, i},\ 2 \leq \ell \leq b, 1 \leq i \leq r;
        \end{cases}
    \end{align}
    and
    \begin{align}
        \xi_{m}(v) = 
        \begin{cases}
            1 - (r-i)\sigma_m, & v = u_{m-1, i},\ 1 \leq i \leq r, \\
            - 1 + (r-i)\sigma_m, & v = u_{m, i},\ 1 \leq i \leq r, \\
            0, & \textotherwise;
        \end{cases}
    \end{align}
    for $m = 3, 4, \ldots, b$. 
    Then $\xi_j$ is a Steklov eigenfunction corresponding to Steklov eigenvalue $\sigma_j$ for $j = 1, 2, \ldots, b$.
\end{proof}

Recall that He-Hua proved an upper bound of the first nontrivial Steklov eigenvalue by the diameter of a tree. 

\begin{theorem}[{\cite[Theorem 1.3]{he2022UpperBoundsSteklov}}]\label{thm:father}
    Let $T = (V, E)$ be a tree with leaves as boundary $B$. 
    Let $D$ be the diameter of $T$. 
    Then 
    \begin{align}
        \sigma_2 \leq \frac{2}{D}.
    \end{align}
    Moreover, if $D$ is odd, then a necessary condition to attain the upper bound is that $T$ is a path. 
    For $i = 1,2, \ldots, D-1$, let $H_i$ be the connected components containing neither $v_0$ nor $v_{D}$ obtained by deleting $v_{i-1}v_i$ and $v_i v_{i+1}$ from edges of $T$. 
    If $D = 2r$ is even, then a necessary condition to attain the upper bound is that $H_i$ is empty for all $1 \leq i \leq D-1$ except for $i = r$. 
\end{theorem}

We characterize when a tree of even diameter achieving the upper bound of first Steklov eigenvalue.  

\begin{lemma}
    Let $T$ be a tree of even diameter $D = 2r$ with leaves as boundary $B$. 
    Suppose $v_0, v_1, \ldots, v_D$ is a diametral path of length $D$. 
    For $i = 1,2, \ldots, D-1$, let $H_i$ be the connected component containing neither $v_0$ nor $v_{D}$ obtained by deleting $v_{i-1}v_i$ and $v_i v_{i+1}$ from edges of $T$. 
    If $H_i$ is empty for all $1 \leq i \leq D-1$ except for $i = r$, then
    $\sigma_2 = 1/r$. 
\end{lemma}

\begin{proof}
    By~\cref{thm:father}, we only need to prove that $\sigma_2 \geq 1/r$. 
    We give the Steklov eigenfunction directly. 
    Let $\sigma_2 = 1/r$ and 
    \begin{align}
        \xi_{2}(v) = 
        \begin{cases}
            1 - i \sigma_2, & v = v_i,\ 0 \leq i \leq 2r, \\
            0, & \textotherwise.
        \end{cases}
    \end{align}
    Then $\xi_2$ is a Steklov eigenfunction corresponding to Steklov eigenvalue $\sigma_2$.
\end{proof}

We use the above graphs to give an lower bound of the value $\sigma_{2, \text{max}}(b, n)$. 

\begin{lemma}\label{lem:party}
    Let $b \geq 2$ and $n \geq b + 1$ be integers. 
    Then the followings hold. 
    \begin{enumerate}
        \item If $b = 2$, then $\sigma_{2, \text{max}}(b, n) = \frac{2}{n-1}$. 
        \item If $b \geq 3$. 
        \begin{enumerate}
            \item If $r$ is a positive integer such that $2r + b - 1 \leq n \leq br + 1$, then $\sigma_{2, \text{max}}(b, n) \geq \frac{1}{r}$. 
            \item If $r$ is a positive integer such that $n = br + 2$, then $\sigma_{2, \text{max}}(b, n) \geq \frac{b}{b(r+1) - 1} = \frac{1}{r + 1 - \frac{1}{b}}$. 
        \end{enumerate}
    \end{enumerate}
\end{lemma}

\begin{proof}
    Suppose $b = 2$. 
    Then the only tree with $2$ leaves and $n$ vertices is the path graph of length $n-1$. 
    Therefore, $\sigma_{2, \text{max}}(2, n) = \sigma_2(P_n)= \frac{2}{n-1}$. 
    In the following we suppose $b \geq 3$. 
    
    We consider a path $P_{2r+1}$ of length $2r$. 
    Then we attach one ends of $(b-2)$ paths of length at most $r$ to the middle point of $P_{2r+1}$. 
    We obtain a tree $T$ with $b$ leaves. 
    Note that the number of vertices of $T$ is at least $(2r + 1 + b - 2 = 2r + b - 1)$ and at most $(b r + 1)$. 
    Let $n$ be the number of vertices of $T$. 
    Then 
    $\sigma_{2, \text{max}}(b, n) \geq \sigma_2(T) = \frac{1}{r}$. 

    We consider the almost fork graph $T = AF(b, r)$. 
    Note that $T$ has $b$ leaves and $n = (br+2)$ vertices. 
    Therefore, $\sigma_{2, \text{max}}(b, n) \geq \sigma_2(T) = \frac{b}{b(r+1) - 1}$.
\end{proof}

We need a graph family larger than the almost fork graph to determine the maximum first Steklov eigenvalue. 

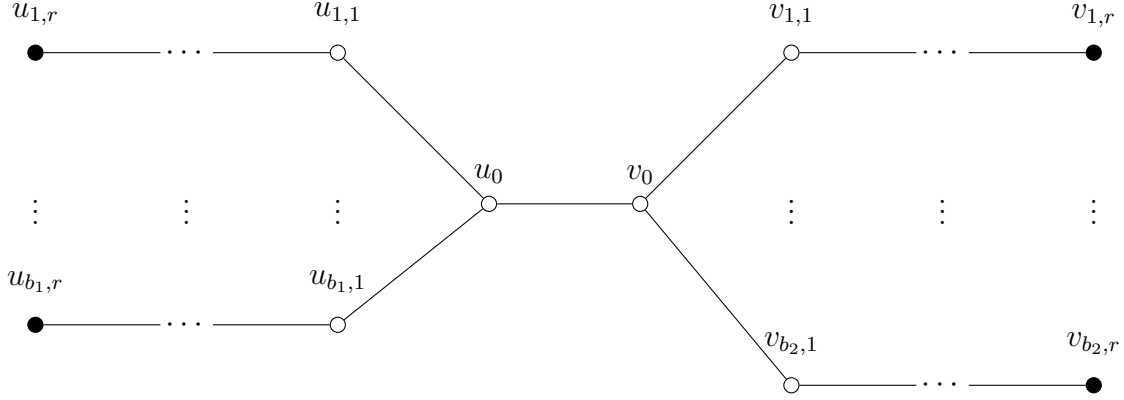
\begin{figure}
    \centering
    \begin{tikzpicture}
    [scale=2, every node/.style={circle, inner sep=1pt, minimum size = 2mm}]

    \node[draw, label=above:{$u_0$}] (u0) at (-1, 0) {};
    \node[draw, label=above:{$u_{1,1}$}] (u11) at (-2, 1) {};
    \node (um1) at (-2, 0) {$\vdots$};
    \node[draw, label=above:{$u_{b_1,1}$}] (ub1) at (-2, -0.8) {};


    \node (u13) at (-3, 1) {$\cdots$};
    \node (um3) at (-3, 0) {$\vdots$};
    \node (ub3) at (-3, -0.8) {$\cdots$};

    \node[draw, fill=black, label=above:{$u_{1,r}$}] (u14) at (-4, 1) {};
    \node (um4) at (-4, 0) {$\vdots$};
    \node[draw, fill=black, label=above:{$u_{b_1,r}$}] (ub4) at (-4, -0.8) {};

    \draw (u0) -- (u11) -- (u13) -- (u14);
    \draw (u0) -- (ub1) -- (ub3) -- (ub4);

    \node[draw, label=above:{$v_0$}] (v0) at (0, 0) {};
    \node[draw, label=above:{$v_{1,1}$}] (v11) at (1, 1) {};
    \node (um1) at (1, 0) {$\vdots$};
    \node[draw, label=above:{$v_{b_2,1}$}] (vb1) at (1, -1.2) {};


    \node (v13) at (2, 1) {$\cdots$};
    \node (vm3) at (2, 0) {$\vdots$};
    \node (vb3) at (2, -1.2) {$\cdots$};

    \node[draw, fill=black, label=above:{$v_{1,r}$}] (v14) at (3, 1) {};
    \node (um4) at (3, 0) {$\vdots$};
    \node[draw, fill=black, label=above:{$v_{b_2,r}$}] (vb4) at (3, -1.2) {};

    \draw (v0) -- (v11) -- (v13) -- (v14);
    \draw (v0) -- (vb1) -- (vb3) -- (vb4);

    \draw (u0) -- (v0);

\end{tikzpicture}
    \caption{a crab graph $CG(b_1, b_2, r)$ with leaves as boundary.}\label{fig:spring}
\end{figure}

\begin{definition}
    Let $b_1, b_2, r \geq 1$ be integers. 
    The crab graph $CG(b_1, b_2, r)$ is obtained by gluing one ends of $b_1$ paths of length $r$ to an end of an edge, and gluing one ends of $b_2$ paths of length $r$ to the other end of the edge.  
    See~\cref{fig:spring}.
    Note that the almost fork graph $AF(b, r) \cong CG(b-1, 1, r) \cong CG(1, b-1, r)$.
\end{definition}

The Steklov eigenvalue of crab graphs are given as follows.

\begin{lemma}
    Let $b_1, b_2, r \geq 1$ be integers. 
    Let $T = (V, E)$ be the crab graph $CG(b_1, b_2, r)$ with leaves as boundary $B$. 
    Then the Steklov eigenvalues of $T$ are
    $\sigma_1 = 0$, $\sigma_2 = \frac{b_1 + b_2}{b_1 b_2 + r(b_1 + b_2)}$, and $\sigma_3 = \cdots = \sigma_{b_1 + b_2} = 1/r$. 
\end{lemma}

\begin{proof}
    We give the Steklov eigenfunctions directly. 
    Define $\sigma_1 = 0$, $\sigma_2 = \frac{b_1 + b_2}{b_1 b_2 + r(b_1 + b_2)}$, and $\sigma_3 = \cdots = \sigma_{b_1 + b_2} = 1/r$. 
    Define
    \begin{align}
        \xi_1(v) = 1, \quad \forall v \in V; 
    \end{align}
    \begin{align}
        \xi_2(v) = 
        \begin{cases}
            b_2(1 - r \sigma_2), & v = u_0, \\
            -b_1(1 - r \sigma_2), & v = v_0, \\
            b_2(1 - (r-i)\sigma_2), & v = u_{\ell, i},\ 1 \leq \ell \leq b_1, 1 \leq i \leq r, \\
            -b_1(1 - (r-i)\sigma_2), & v = v_{\ell, i},\ 1 \leq \ell \leq b_2, 1 \leq i \leq r;
        \end{cases}
    \end{align}
    \begin{align}
        \xi_{m}(v) = 
        \begin{cases}
            1 - (r-i)\sigma_m, & v = u_{m-2, i},\ 1 \leq i \leq r, \\
            - 1 + (r-i)\sigma_m, & v = u_{m - 1, i},\ 1 \leq i \leq r, \\
            0, & \textotherwise;
        \end{cases}
    \end{align}
    for $m = 3, 4, \ldots, b_1 + 1$, and
    \begin{align}
        \xi_{m}(v) = 
        \begin{cases}
            1 - (r-i)\sigma_m, & v = v_{m - b_1 - 1, i},\ 1 \leq i \leq r, \\
            - 1 + (r-i)\sigma_m, & v = v_{m - b_1, i},\ 1 \leq i \leq r, \\
            0, & \textotherwise
        \end{cases}
    \end{align}
    for $m = b_1 + 2, b_1 + 3, \ldots, b_1 + b_2$.
    Then $\xi_j$ is a Steklov eigenfunction corresponding to Steklov eigenvalue $\sigma_j$ for $j = 1, 2, \ldots, b$.
\end{proof}

Now we are prepared to determine the maximum first Steklov eigenvalue. 

\begin{proof}[{Proof of~\cref{thm:question}}]
    Note that $\cT(2, n) = \set{P_{n}}$, therefore~\cref{eq:parent} holds. 
    For $b \geq 3$ and $n \geq b+1$. 
    Let $T \in \cT(b, n)$ be a tree with $b$ leaves and $n$ vertices. 

    Suppose $b(r-1) + 3 \leq n \leq br+1$ for some positive integer $r \in \ZZ_+$. 
    Since $b(r-1) + 3 \geq 2r + b - 1$, by~\cref{lem:party} we have $\sigma_{2, \text{max}}(b, n) \geq \frac{1}{r}$. 
    We claim that the diameter of $T$ is at least $2r$. 
    Then $\sigma_{2, \text{max}}(b, n) \leq \frac{2}{D} \leq \frac{1}{r}$, which implies that $\sigma_{2, \text{max}}(b, n) = \frac{1}{r}$. 
    In fact, if the diameter $D$ is at most $(2r-1)$, then the number of vertices $n$ is at most $b(r-1) + 1 + 1 = b(r-1) + 2$. 
    Contradiction. 

    Suppose $n = br + 2$ for some positive integer $r \in \ZZ_+$. 
    By~\cref{lem:party} we have $\sigma_{2, \text{max}}(b, n) \geq \frac{1}{r + 1 - \frac{1}{b}}$. 
    We claim that the diameter of $T$ is at least $2r + 1$. 
    In fact, if the diameter $D$ is at most $2r$, then the number of vertices $n$ is at most $br + 1$. 
    Contradiction. 
    Suppose the diameter of $T$ is at least $2r+2$. 
    Then $\sigma_{2}(T) \leq \frac{2}{2r+2} \leq \frac{1}{r+1} < \frac{1}{r+1 - \frac{1}{b}}$.
    We further claim that the tree of diameter $2r+1$ with $b$ leaves and $br + 2$ vertices are exactly the crab graphs $CG(b-1, 1, r), CG(b-2, 2, r), \ldots, CG(1, b-1, r)$. 
    Note that $\sigma_2(CG(b-\ell, \ell, r)) = \frac{1}{r + \frac{\ell(b-\ell)}{b}} \leq \frac{1}{r+1 - \frac{1}{b}}$ for $\ell = 1,2, \ldots, b-1$ implies that $\sigma_2(b, br+2) = \frac{1}{r+1 - \frac{1}{b}}$. 
    In fact, consider a diametral path of length $2r+1$. 
    For the rest of $(b-2)$ leaves, consider their distance to this diametral path. 
    This distance $d$ is at most $r$, which would contribute $d (\leq r)$ vertices. 
    Since the total number of vertices is $2r+2 + (b-2)r = br+2$, the only possibility is that the distance is exactly $r$ for each of the rest $(b-2)$ leaves. 
    In such case, these leaves can only be connected to the middle two vertices of the diametral path. 
    So the graph is a crab graph $CG(b-\ell, \ell, r)$ for some $\ell = 1,2, \ldots, b-1$.
\end{proof}

For $k \geq 3$, the maximum $k$-th Steklov eigenvalue behaves quite differently. 

\begin{definition}\label{def:barbell}
    Attach $p$ pendents to an end of a path of length $D-2$ and $q$ pendents to the other end of the path. 
    We get a graph called the barbell graph $B(p,q,D)$,
    (see \cref{fig:barbell}).
    \begin{figure}[htbp]
        \centering
        \begin{tikzpicture}
    [scale=2, every node/.style={circle, inner sep=1pt, minimum size = 2mm}]
    \node[draw, fill=black, label=left:{$u_1$}] (u1) at (-1, 1) {};
    \node[draw, fill=black, label=left:{$u_2$}] (u2) at (-1, 0.8) {};
    \node (u2-a) at (-1,0.2) {$\vdots$};
    \node[draw, fill=black, label=left:{$u_p$}] (ua) at (-1, -0.4) {};
    
    \node[draw, label=above:{$w_1$}] (w1) at (0, 0) {};
    \node[draw, label=above:{$w_2$}] (w2) at (1, 0) {};
    \node[draw, label=above:{$w_3$}] (w3) at (2, 0) {};
    \node (w2-3) at (2.5,0) {$\cdots$};
    \node[draw, label=above:{$w_{L-3}$}] (w4) at (3, 0) {};
    \node[draw, label=above:{$w_{L-2}$}] (w5) at (4, 0) {};
    \node[draw, label=above:{$w_{L-1}$}] (w6) at (5, 0) {};
    
    \node[draw, fill=black, label=right:{$v_1$}] (v1) at (6, 1) {};
    \node[draw, fill=black, label=right:{$v_2$}] (v2) at (6, 0.8) {};
    \node (v2-b) at (6,-0.1) {$\vdots$};
    \node[draw, fill=black, label=right:{$v_q$}] (vb) at (6, -1) {};

    \draw (u1) -- (w1);
    \draw (u2) -- (w1);
    \draw (u2-a) -- (w1);
    \draw (ua) -- (w1);
    
    \draw (w1) -- (w2) -- (w3);
    \draw (w4) -- (w5) -- (w6);
    
    \draw (v1) -- (w6);
    \draw (v2) -- (w6);
    \draw (v2-b) -- (w6);
    \draw (vb) -- (w6);

\end{tikzpicture}
        \caption{The barbell graph $B(p,q,D)$ with leaves as boundary.}
        \label{fig:barbell}
    \end{figure}
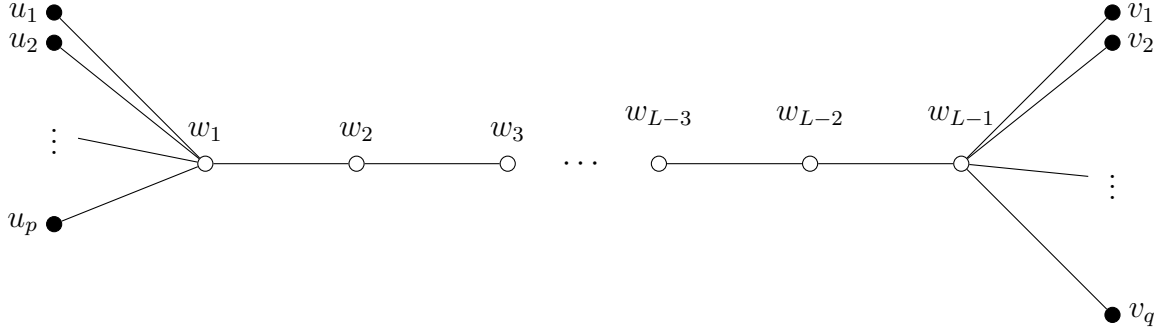
\end{definition}

Recall the Steklov eigenvalues of the barbell graphs. 

\begin{lemma}[{\cite[Lemma 2.10]{linFirstSteklovEigenvalue2024} or~\cite[Lemma 2.4]{yu2022MinimalSteklovEigenvalues}}]\label{lem:barbell}
    Let $T$ be the barbell graph $B(p,q,D)$ with leaves as boundary. 
    Then $\sigma_1 = 0$, $\sigma_2= \frac{p+q}{(D-2)p q + p + q}$, and $\sigma_3 = \cdots = \sigma_{p+q} = 1$.
\end{lemma}

For completeness, we repeat the proof of~\cref{lem:barbell} here. 

\begin{proof}
    Label the vertices as in \cref{fig:barbell}. 
    We give a complete set of harmonic extensions of eigenfunctions for $T$. 
    Define 
    \begin{align}
        f_1(x) = 1, x \in V.
    \end{align}
    Then $f_1$ is the harmonic extensions of eigenfunction corresponding to eigenvalue $0$.
    Define 
    \begin{align}
        f_2(x) = 
        \begin{cases}
            1, & x = u_i, i = 1,2, \ldots, a, \\
            1 - \dfrac{p+q}{(D-2)p q + p + q} + (i-1)\dfrac{2pq}{(D-2)pq + p + q}, & x = w_i, i = 1,2, \ldots, L-2,\\
            -1, & x = v_i, i = 1,2, \ldots, b.
        \end{cases}
    \end{align}
    Then $f_2$ is the harmonic extensions of eigenfunction corresponding to eigenvalue $\dfrac{p+q}{(D-2)p q + p + q}$. 
    Define 
    \begin{align}
        f_m(x) = 
        \begin{cases}
            1, & x = u_1, \\
            -1, & x = u_m, \\
            0, & \textotherwise.
        \end{cases}
    \end{align}
    Then $f_m$ is the harmonic extensions of eigenfunction corresponding to eigenvalue $1$ for $m = 2,3, \ldots, p$.  
    Define 
    \begin{align}
        f_{p+m}(x) = 
        \begin{cases}
            1, & x = v_1, \\
            -1, & x = v_m, \\
            0, & \textotherwise.
        \end{cases}
    \end{align}
    Then $f_{p+m}$ is the harmonic extensions of eigenfunction corresponding to eigenvalue $1$ for $m = 2,3, \ldots, q$. 
    Since there are $p+q$ eigenvalues in total, we know that the eigenvalues are $\displaystyle 0, \frac{p+q}{(D-2)p q + p + q}, 1, \ldots, 1$. 
    And hence $\displaystyle \sigma_2 = \frac{p+q}{(D-2)p q + p + q}$.
\end{proof}

The following theorem provides uniform upper bound for the Steklov eigenvalues. 

\begin{theorem}[{\cite[Corollary 1.6]{lin2024UpperBoundsSteklov}}]\label{thm:degree_bound}
    Let $G = (V,E)$ be a graph with boundary $B$ such that the boundary vertices form an independent set of $G$. 
    Suppose that the degree sequence of the boundary vertices is $d_1 \leq d_2 \leq \cdots \leq d_{\card{B}}$. 
    Then
    \begin{align*}
        \sigma_k \leq d_k,\quad k = 1,2, \ldots, \card{B}. 
    \end{align*}
\end{theorem}

For completeness, we repeat the proof of~\cref{thm:degree_bound} here. 

Let $L$ be the Laplacian matrix of the graph $G = (V, E)$. 
Let $B$ be the boundary of the graph $G$. 
For every positive real number $r$, we define a diagonal matrix $D^{(r)}$ by
\begin{align*}
    D^{(r)}(x,x) = 
    \begin{cases}
        1, & x \in B, \\
        r, & \textotherwise.
    \end{cases}
\end{align*}
Next we consider the eigenvalues of the matrix $L^{(r)} = D^{(r)} L D^{(r)}$.

\begin{lemma}[{\cite[Proposition 3]{hassannezhad_higher_2020}}]\label{lem:steklov_limit}
    Let $G = (V, E)$ be a graph with boundary $B$. 
    For every positive real number $r$, the matrix $L^{(r)}$ is defined as above.
    Let $\mu_1^{(r)} \leq \mu_2^{(r)} \leq \cdots \leq \mu_{|V|}^{(r)}$ be the eigenvalues of $L^{(r)}$. 
    Then
    \begin{align*}
        \lim_{r \to +\infty} \mu_k^{(r)} = \sigma_k, \quad 1 \leq k \leq \card{B},
    \end{align*}
    and
    \begin{align*}
        \lim_{r \to +\infty} \mu_k^{(r)} = +\infty, \quad \card{B} < k \leq \card{V}.
    \end{align*}
\end{lemma}


We are prepared to prove~\cref{thm:degree_bound}.

\begin{proof}[{Proof of \cref{thm:degree_bound}}]
    Take $M = L^{(r)}$. 
    Let $N$ be the principal submatrix of $M$ indexed by $B$. 
    By Cauchy interlacing theorem, the eigenvalues of $N$ interlace the eigenvalues of $M$. 
    Therefore
    \begin{align*}
        \mu_k(L^{(r)}) \leq \mu_k(N) = d_k, \quad k = 1,2, \ldots, \card{B}.
    \end{align*}
    The conclusion follows by~\cref{lem:steklov_limit}.
\end{proof}

\begin{corollary}\label{coro:tree_steklov_upper}
    Let $T = (V, E)$ be a tree with leaves as boundary $B$. 
    Then $0 = \sigma_1 < \sigma_2 \leq \cdots \leq \sigma_{\card{B}} \leq 1$. 
\end{corollary}

Now we are prepared to determine $\sigma_{k, \text{max}}(b, n)$ for $k \geq 3$. 

\begin{proof}[{Proof of~\cref{thm:natural}}]
    Consider the barbell graph $B(b-1, 1, n-b + 1)$. 
    Then by~\cref{lem:barbell} we have $\sigma_{k, \text{max}}(b, n) \geq \sigma_k(B(b-1, 1, n-b + 1)) = 1$ for $k \geq 3$. 
    On the other hand, by~\cref{coro:tree_steklov_upper} we have $\sigma_{k, \text{max}}(b, n) \leq 1$. 
    Therefore, $\sigma_{k, \text{max}}(b, n) = 1$ for $k \geq 3$.
\end{proof}

On the other hand, if we fix the diameter and the total number vertices of the graph, the maximum first Steklov eigenvalue is not easy to determine. 

\begin{figure}[htbp]
    \begin{tikzpicture}
    [scale=2, every node/.style={circle, inner sep=1pt, minimum size = 2mm}]

    \node[draw, label=above:{$o_0$}] (o0) at (0, 0) {};
    
    \node[draw, label=above:{$u_1$}] (u1) at (-1, 0) {};
    \node (um) at (-2, 0) {$\cdots$};
    \node[draw, fill=black, label=above:{$u_r$}] (ur) at (-3, 0) {};

    \node[draw, label=above:{$v_1$}] (v1) at (1, 0) {};
    \node[draw, label=above:{$v_2$}] (v2) at (2, 0) {};
    \node (vm) at (3, 0) {$\cdots$};
    \node[draw, fill=black, label=above:{$v_{r+1}$}] (vr) at (4, 0) {};

    \node[draw, label=left:{$w_{1,1}$}] (w11) at (-1, -1) {};
    \node (w1m) at (-1, -2) {$\vdots$};
    \node[draw, fill=black, label=left:{$w_{1,c}$}] (w1c) at (-1, -3) {};

    \node[draw, label=right:{$w_{b-2,1}$}] (wb1) at (1, -1) {};
    \node (wbm) at (1, -2) {$\vdots$};
    \node[draw, fill=black, label=right:{$w_{b-2,c}$}] (wbc) at (1, -3) {};

    \node (wm1) at (0, -1) {$\cdots$};
    \node (wmm) at (0, -2) {$\cdots$};
    \node (wmc) at (0, -3) {$\cdots$};

    \draw (o0) -- (u1) -- (um) -- (ur);
    \draw (o0) -- (v1) -- (v2) -- (vm) -- (vr);
    \draw (o0) -- (w11) -- (w1m) -- (w1c);
    \draw (o0) -- (wb1) -- (wbm) -- (wbc);
\end{tikzpicture}
    \caption{an almost seesaw graph $AS(r, b, c)$ with leaves as boundary}\label{fig:almost_seesaw}
\end{figure}
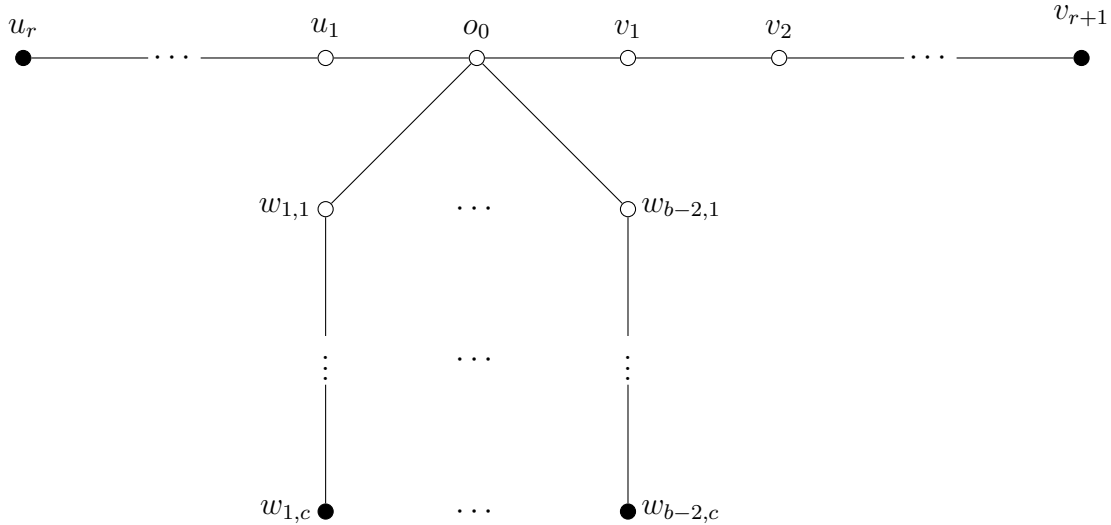

\begin{definition}
    Let $r \geq 1$, $b \geq 2$, $c \geq 0$ be integers. 
    The almost seesaw graph $AS(r, b, c)$ is obtained by gluing one ends of $b$ paths together, where one path is of length $r$, one path is of length $r+1$, and the rest $(b-2)$ paths are of length $c$. 
    See~\cref{fig:almost_seesaw}.
\end{definition}

We give the Steklov eigenvalues of the almost seesaw graphs. 

\begin{lemma}\label{lem:almost_seesaw}
    Let $r \geq 1$, $b \geq 2$, $0 \leq c \leq r$ be integers. 
    Let $T = (V, E)$ be the almost seesaw graph $AS(r, b, c)$ with leaves as boundary $B$. 
    Then the Steklov eigenvalues of $T$ are
    $\sigma_1 = 0$, $\sigma_2 = C^-(r, b, c)$, $\sigma_3 = C^+(r, b, c)$, and $\sigma_4 = \cdots = \sigma_{b} = 1/c$, where 
    \begin{align}
        C^\pm(r, b, c) = \frac{2br + b + 2c - 2r - 1 \pm \sqrt{b^2-2 b+4 c^2-8 c r-4 c+4 r^2+4 r+1}}{2 \left(b r^2+b r+2 c r+c-2 r^2-2 r\right)}.
    \end{align}
\end{lemma}

\begin{proof}
    We give the Steklov eigenfunctions directly. 
    Define $\sigma_1 = 0$, $\sigma_2 = C^-(r, b, c)$, $\sigma_3 = C^+(r, b, c)$, and $\sigma_4 = \cdots = \sigma_{b} = 1/c$, and
    \begin{align}
        \alpha^\pm &= \frac{1}{2} \left(\pm\sqrt{b^2-2 b+4 c^2-8 c r-4 c+4 r^2+4 r+1}-b-2 c+2 r+3\right), \\
        \beta^\pm &= \frac{1}{2} \left(\pm\sqrt{b^2-2 b+4 c^2-8 c r-4 c+4 r^2+4 r+1}-b+2 c-2 r+1\right)
    \end{align}
    Define
    \begin{align}
        \xi_1(v) = 1, \quad \forall v \in V; 
    \end{align}
    \begin{align}
        \xi_2(v) = 
        \begin{cases}
            \alpha^+(1 - (r-i)\sigma_2), & v = u_i,\ 1 \leq i \leq r, \\
            \beta^-(1 - (r+1-i)\sigma_2), & v = v_i,\ 1 \leq i \leq r, \\
            1 - (c-i)\sigma_2, & v = w_{\ell, i},\ 1 \leq \ell \leq b-2, 1 \leq i \leq c, \\
            1 - c \sigma_2, & v = o_0;
        \end{cases}
    \end{align}
    \begin{align}
        \xi_3(v) = 
        \begin{cases}
            \alpha^-(1 - (r-i)\sigma_3), & v = u_i,\ 1 \leq i \leq r, \\
            \beta^+(1 - (r+1-i)\sigma_3), & v = v_i,\ 1 \leq i \leq r, \\
            1 - (c-i)\sigma_3, & v = w_{\ell, i},\ 1 \leq \ell \leq b-2, 1 \leq i \leq c, \\
            1 - c \sigma_3, & v = o_0;
        \end{cases}
    \end{align}
    and
    \begin{align}
        \xi_{m}(v) = 
        \begin{cases}
            1 - (c-i)\sigma_m, & v = w_{m - 3, i},\ 1 \leq i \leq r, \\
            - 1 + (c-i)\sigma_m, & v = w_{m - 2, i},\ 1 \leq i \leq r, \\
            0, & \textotherwise
        \end{cases}
    \end{align}
    for $m = 4, 5, \ldots, b$.
    Then $\xi_j$ is a Steklov eigenfunction corresponding to Steklov eigenvalue $\sigma_j$ for $j = 1, 2, \ldots, b$.
\end{proof}

Recall that shrinking a tree only enlarges the Steklov eigenvalues. 

\begin{lemma}[{\cite[Corollary 1.1]{yu2024MonotonicitySteklovEigenvalues}}]\label{lem:monotone}
    Let $T$ be a finite tree with leaves as boundary $B$. 
    Let $T'$ be a subtree of $T$ with leaves as boundary $B'$. 
    Then
    \begin{align}
        \sigma_{i}(T) \leq \sigma_{i}(T')
    \end{align}
    for $i = 1,2, \ldots, \card{B'}$.
\end{lemma}

The maximum value of $\widetilde{\sigma}_2(3, n)$ is given by the barbell graph $B(1, n-3, 3)$.  

\begin{lemma}\label{lem:diameter_3}
    Let $n \geq 4$ be an integer. 
    Then
    \begin{align}
        \widetilde{\sigma}_{2, \text{max}}(3, n) = \frac{n-2}{2n - 5}.
    \end{align}
    The tree attaining the bound is the barbell graph $B(1, n-3, 3)$. 
\end{lemma}

\begin{proof}
    Note that all trees of diameter three is a barbell graph $B(p, q, 3)$. 
    If the number of vertices is $n$, then $p + q = n - 2$. 
    \begin{align}
        \widetilde{\sigma}_{2, \text{max}}(3, n) = \max_{p + q = n-2} \sigma_2(B(p, q, 3)) = \max_{p+q = n-2} \frac{p+q}{pq + p + q} = \frac{n-2}{2n-5}. 
    \end{align}
\end{proof}

We use almost seesaw graphs to give lower bounds for $\widetilde{\sigma}_{2, \text{max}}$. 

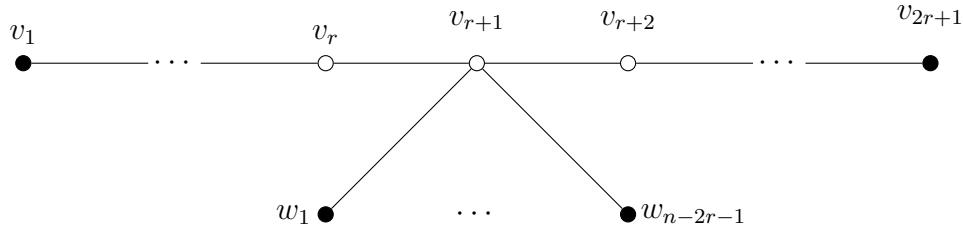
\begin{figure}[htbp]
    \centering
    \begin{tikzpicture}
    [scale=2, every node/.style={circle, inner sep=1pt, minimum size = 2mm}]

    \node[draw, label=above:{$v_{r+1}$}] (o0) at (0, 0) {};
    
    \node[draw, label=above:{$v_{r}$}] (u1) at (-1, 0) {};
    \node (um) at (-2, 0) {$\cdots$};
    \node[draw, fill=black, label=above:{$v_1$}] (ur) at (-3, 0) {};

    \node[draw, label=above:{$v_{r+2}$}] (v1) at (1, 0) {};
    \node (vm) at (2, 0) {$\cdots$};
    \node[draw, fill=black, label=above:{$v_{2r+1}$}] (vr) at (3, 0) {};

    \node[draw, fill=black, label=left:{$w_{1}$}] (w11) at (-1, -1) {};

    \node[draw, fill=black, label=right:{$w_{n-2r-1}$}] (wb1) at (1, -1) {};

    \node (wm1) at (0, -1) {$\cdots$};

    \draw (o0) -- (u1) -- (um) -- (ur);
    \draw (o0) -- (v1) -- (vm) -- (vr);
    \draw (o0) -- (w11);
    \draw (o0) -- (wb1);
\end{tikzpicture}
    \caption{a tree of even diameter $D=2r$ with $\sigma_2 = 2/D$.}\label{fig:ruler}
\end{figure}

\begin{proof}[{Proof of~\cref{thm:gas}}]
    Note that $\widetilde{\cT}(1, 2) = \set{P_{2}}$ and $\widetilde{\cT}(1, n) = \emptyset$ for $n \geq 3$, therefore~\cref{eq:order} holds. 
    For $D \geq 2$ and $n \geq b+1$. 
    Let $T \in \widetilde{\cT}(D, n)$ be a tree of diameter $D$ with $n$ vertices. 

    Suppose $D = 2r$ is even. 
    Consider a path $P_{2r+1}$ of length $2r$. 
    Then we attach one ends of $(n-2r-1)$ edges to the middle point of $P_{2r+1}$. 
    We obtain a tree $T$ of diameter $D$ with $n$ vertices, see~\cref{fig:ruler}. 
    Then $\widetilde{\sigma}_{2, \text{max}}(D, n) \geq \sigma_2(T) = 2/D$. 
    On the other hand, by~\cref{thm:father} we have $\widetilde{\sigma}_{2, \text{max}}(D, n) \leq 2/D$. 
    Therefore, $\widetilde{\sigma}_{2, \text{max}}(D, n) = 2/D$.

    For $D = 3$, we conclude by~\cref{lem:diameter_3}.

    Suppose $D = 2r+1$ is odd. 
    If $2r+2 \leq n \leq 2r + (b-2)c + 2$. 
    Consider the almost seesaw graph $AS(r, b, c)$. 
    There is a subtree $T$ of $AS(r, b, c)$ of diameter $D$ with $n$ vertices. 
    Then by~\cref{lem:almost_seesaw,lem:monotone} we have $\widetilde{\sigma}_{2, \text{max}}(D, n) \geq \sigma_2(T) \geq \sigma_2(AS(r,b,c)) = C^-(r, b, c)$. 
\end{proof}

\begin{conjecture}
    Let $D = 2r+1$ be an odd integer and $n \geq D+1$ an integer. 
    Then
    \begin{align}
        \widetilde{\sigma}_{2, \text{max}}(D, n) = C^-(r, n - 2r - 2, 1).
    \end{align}
    The tree attaining the bound is the almost seesaw graph $AS(r, n - 2r - 2, 1)$. 
\end{conjecture}

\section*{Acknowledgements}
Huiqiu Lin was supported by the National Natural Science Foundation of China (No. 12271162, No. 12326372), and the Natural Science Foundation of Shanghai (No. 22ZR1416300 and No. 23JC1401500) and The Program for Professor of Special Appointment (Eastern Scholar) at Shanghai Institutions of Higher Learning (No. TP2022031). 
Da Zhao was supported in part by the National Natural Science Foundation of China (No. 12471324), and the Natural Science Foundation of Shanghai, Shanghai Sailing Program (No. 24YF2709000).

\bibliographystyle{plain}
\bibliography{ref}

\end{document}